\documentclass[12pt]{amsart}

\usepackage{amssymb, amsmath, amsthm, hyperref}

\numberwithin{equation}{section}

\newtheorem{prop}{Proposition:}[section]
\newtheorem{cor}[prop]{Corollary:}

\newtheorem{thm}[prop]{Theorem:}
\newtheorem{lemma}[prop]{Lemma:}

\theoremstyle{definition}

\newcommand{\del}{\partial}

\newcommand{\co}[1]{\til{#1}^r}
\newcommand{\cn}[1]{\til{\underline{#1}}}
\newcommand{\delr}[1]{\frac{\del^{#1}}{\del r^{#1}}}
\newcommand{\bt}{\boxtimes}
\newcommand{\nm}[2]{|| #1 ||_{#2}}
\newcommand{\nmm}[3]{\brs{#1}_{#2;#3}}
\newcommand{\TT}{\mathcal T}
\newcommand{\TTs}[3]{\TT^{#1}_{#2;#3}}
\renewcommand{\ge}{\epsilon}

\newcommand{\til}[1]{\widetilde{#1}}

\newcommand{\gt}{\theta}
\newcommand{\ga}{\alpha}
\newcommand{\gs}{\sigma}
\newcommand{\gl}{\lambda}
\newcommand{\gb}{\beta}
\newcommand{\gd}{\delta}
\newcommand{\gk}{\kappa}
\newcommand{\gw}{\omega}
\newcommand{\gp}{\phi}

\newcommand{\gD}{\Delta}

\newcommand{\WW}{\mathcal W}

\newcommand{\brs}[1]{\left|#1 \right|}
\newcommand{\N}{\nabla}
\newcommand{\up}{\uparrow}
\newcommand{\dn}{\downarrow}
\newcommand{\PP}{\mathcal P}
\newcommand{\OO}{\mathcal O}

\DeclareMathOperator{\Rm}{Rm}
\DeclareMathOperator{\W}{W}
\DeclareMathOperator{\Rc}{Rc}

\DeclareMathOperator{\Vol}{Vol}

\DeclareMathOperator{\sym}{sym}
\DeclareMathOperator{\tr}{tr}

\begin{document}
\title[Removable Singularities of Bach-Flat Metrics]
{Asymptotic Curvature Decay and Removal of Singularities of Bach-Flat
Metrics}
\author{Jeffrey Streets}
\address{Fine Hall\\
         Princeton University\\
         Princeton, NJ 08544}
\email{\href{mailto:jstreets@math.princeton.edu}{jstreets@math.princeton.edu}}

\date{August 7, 2007}

\begin{abstract}
\noindent We prove a removal of singularities result for Bach-flat
metrics in dimension $4$ under the assumption of bounded $L^2$ norm
of curvature, bounded Sobolev constant and a volume growth bound.
This result extends the removal of singularities result for special
classes of Bach-flat metrics obtained in \cite{TVMOD}.  For the
proof we analyze the decay rates of solutions to the Bach-flat
equation linearized around a flat metric. This classification is
used to prove that Bach-flat cones are in fact ALE of order $\tau$
for any $\tau < 2$.  This result is then used to prove the removal
of singularities theorem.
\end{abstract}

\maketitle

\section{Introduction}
In dimension $4$, the Euler-Lagrange equations for the functional
\begin{gather} \label{Wdef}
\WW(g) := \int_M \brs{\W_g}^2 dV_g
\end{gather}
where $\W_g$ is the Weyl curvature, are given by
\begin{gather} \label{Bdef}
B_{ij} = \N^k \N^l W_{ikjl} + \frac{1}{2} R^{kl} W_{ikjl} = 0
\end{gather}
where $W_{ijkl}$ and $R^{kl}$ correspond to the components of the Weyl and
Ricci tensors respectively \cite{Derd}, \cite{Besse}.  $B_{ij}$ are the
components of the
\emph{Bach tensor}, and a metric which is critical for $\WW$ is called
\emph{Bach-flat}.

A smooth Riemannian manifold $(X^4, g)$ is called an asymptotically locally
Euclidean end of order $\tau$ if there exists a finite subgroup $\Gamma \subset
SO(4)$  which acts freely on $\mathbb R^4 \backslash B(0, R)$ and a $C^\infty$
diffeomorphism $\phi : X \to \left( \mathbb R^4 \backslash B(0, R) \right) /
\Gamma$ such that using this identification
\begin{align}
g_{ij} =&\ \gd_{ij} + O(r^{-\tau})\\
\del^{|k|} g_{ij} =&\ O(r^{- \tau - k})
\end{align}
for any partial derivative of order $k$ as $r \to \infty$.  We say an end is ALE
of order $0$ if there exist coordinates as above so that
\begin{align}
g_{ij} =&\ \gd_{ij} + o(1)\\
\del^{|k|} g_{ij} =&\ o(r^{-k})
\end{align}
as $r \to \infty$.

\begin{thm} \label{decaythm} Let $(X^4, g)$ be a complete, noncompact
four-dimensional Bach-flat space with zero scalar curvature which is ALE of
order $0$.  Assume further that
\begin{gather}
\int_X \brs{\Rm_g}^2 dV_g < \Lambda < \infty, C_S < K < \infty
\end{gather}
Then the cone is ALE of order $\tau$ for any $\tau < 2$.
\end{thm}
This theorem is an extension of theorem 1.1 of \cite{TVBF}, where under the same
hypotheses on curvature and Sobolev constant and also assuming $b_1(X) < \infty$
where $b_1(X)$ is the first Betti number it is proved that the end is ALE of
order zero.

\begin{thm} \label{maintheorem} Let $(X, d, x)$ be a complete locally compact
length space with base point $x$.  Let $B(x,1) \setminus \{x \}$ be a
connected $C^\infty$ Bach-flat four-manifold satisfying
\begin{align}
& \int_{B(x,1)} \brs{\Rm_g}^2 dV_g < \infty\\
& \nm{u}{L^4(B(x,1) \setminus \{x \})} \leq C_s \nm{\N u}{L^2(B(x,1) \setminus
\{x\})}, u \in C^{0,1}(B(x,1) \setminus \{x \})\\
& \Vol(B(x,r)) \leq V_1 r^4, r > 0\\
& b_1(X) < \infty
\end{align}
where $C_s$ and $V_1$ are positive constants.  Then the metric $g$ extends to
$B(x,1)$ as a smooth orbifold metric.
\end{thm}

This is an extension of theorem 6.4 of \cite{TVMOD}, and indeed this result was
suggested in the work of Gang Tian and Jeff Viaclovsky \cite{TVMOD},
\cite{TVBF}.  There removal of singularities is
proved for certain subclasses of Bach-flat manifolds, specifically
half-conformally flat metrics, metrics with harmonic curvature, and constant
scalar curvature K\"ahler metrics.  The proof uses certain Kato inequalities
which are satisfied in these special cases to improve the decay of the Ricci
tensor, which in turn can be used to get improved decay of the full curvature
tensor.  The techniques used herin also suffice to give this result assuming
that the
metric has harmonic curvature, that is,
\begin{gather} \label{harmoniccurvature}
\gd \Rm = - \N^i R_{ijkl} = 0
\end{gather}

We now give an outline of the rest of the paper.  In section $2$ we decompose
the action of the Bilaplacian acting on symmetric
two-tensors according to the radial separation of variables on a cone.  In
section $3$ we compute the linearization of the Bach-flat condition around a
flat metric.  We observe that after applying a change of diffeomorphism gauge
and possibly a conformal change this equation is equivalent to the equation
$\gD^2 h = 0$ where $h$ is the perturbation of the metric.  Using the separation
of variables we further simplify this equation and classify the solutions in
terms of their decay or growth rates at infinity.  Section $4$ gives the proof
of theorem \ref{decaythm}, which consists of considering the Bach-flat equation
at infinity in the cone as a perturbation of a flat metric, and using our
analysis of the linearized equation.  Section $5$ uses theorem \ref{decaythm} to
prove theorem \ref{maintheorem}.

\textbf{Acknowledgements:} The author would like to express deep gratitude to
Gang Tian for suggesting this problem and suggesting that the arguments of
\cite{Tian} could be used to prove theorem \ref{maintheorem}.

\section{Decomposition of the Bilaplacian}
In this section we will decompose the action of the Bilaplacian on symmetric
two-tensors on $\mathbb R^4$ according to the cone structure.  These tensors
naturally decompose according to radial and tangential directions.  We
reintroduce the notation of \cite{Tian} and recall the decomposition of the
Laplacian proved there.  These lemmas are then used in a straightforward way to
decompose the action of the Bilaplacian.  The calculation is greatly simplified
by restricting to $\mathbb R^4$, where we completely understand the curvature
tensor, and how it decomposes according to the cone structure.

Write $\left(\mathbb R^4, g_{\mbox{flat}} \right)$ as $C(S^3, g_{\mbox{can}})$,
the cone over the three-sphere with its canonical metric.  Let $P: T(1, S^3) \to
T(r, S^3)$ denote the
identification of tangent bundles induced by parallel transport along radial
geodesics, and let $\til{\N}^r$ denote the Riemannian connection with respect to
the induced metric on $\left(r, S^3 \right)$.  Define the connection
\begin{gather} \label{conndef}
\cn{\N} = P \til{\N}^1 P^{-1}
\end{gather}
A basic formula is then
\begin{gather}
\co{\N} = r^{-1} \cn{\N}
\end{gather}
Also let $\co{d}, \co{\gd}, \co{\tr}$ be defined with respect to the induced
metric on $(r, S^3)$ and let $\co{d}$ denote exterior differentiation on
$(r, S^3)$.  Then define
\begin{gather}
\begin{split}
\cn{\gd} =&\ P \til{\gd}^1 P^{-1}\\
\cn{d} =&\ P \til{d}^1 P^{-1}\\
\cn{\tr} =&\ P \til{\tr}^1 P^{-1}\\
\end{split}
\end{gather}
The following formulae are immediate
\begin{gather}
\begin{split}
\co{\gD} =&\ r^{-2} \cn{\gD}\\
\co{\gd} =&\ r^{-1} \cn{\gd}\\
\co{d} =&\ r^{-1} \cn{d}\\
\co{\tr} =&\ \cn{\tr}
\end{split}
\end{gather}
Now say $e$ is a vector tangent to $(r, S^3)$ and $e^*$ is its dual $1$-form.
Let $\eta$ be a $1$-form such that
\begin{gather} \label{decloc20}
\eta \left(\frac{\del}{\del r} \right) = 0\\ \label{decloc30}
\N_{\frac{\del}{\del r}} \eta \equiv 0
\end{gather}
Some simple formulas are then
\begin{align}
\N_e \eta =&\ r^{-1} \left( \cn{\N}_e \eta - \eta(e) dr \right)\\
\N_e dr =&\ r^{-1} e^*\\
\N_{\frac{\del}{\del r}} \frac{\del}{\del r} =&\ \N_{\frac{\del}{\del r}} dr =
0
\end{align}
Now fix $\delr{}, e_1, e_2, e_3$ a local orthonormal basis satisfying
\begin{gather} \label{decloc40}
\co{\N}_{e_i} e_j = 0
\end{gather}
at a fixed point $(r, x)$ and also
\begin{gather}
\N_{\delr{}} e \equiv 0
\end{gather}
Finally, given $\gw_1, \gw_2$ $1$-forms let
\begin{gather}
\gw_1 \boxtimes \gw_2 = \gw_1 \otimes \gw_2 + \gw_2 \otimes \gw_1
\end{gather}

Using this frame the symmetric two-tensors naturally decompose into three
distinct types.  In the next four lemmas we record the action of divergence and
the Laplacian on these tensors.
We will need these preliminary calculations
to calculate the action of the Bilaplacian on each of these types of tensors.
The proofs of these lemmas can all be found in \cite{Tian}.

\begin{lemma} Consider $f(r) B$ where $B = \sum \eta_i \bt \eta_j$ and $\eta_i$
satisfies (\ref{decloc20}) and (\ref{decloc30}).  Then
\begin{gather} \label{horizlap}
\begin{split}
\N^* \N (fB) =&\ \left( - f'' - 3 r^{-1} f' + r^{-2} f \left(\cn{\N}^* \cn{\N}
+ 2 \right)\right) B\\
&\ + 2 r^{-2} f \cn{\gd}(B) \bt dr - 2 r^{-2} f \cn{\tr}(B) dr \otimes dr
\end{split}
\end{gather}
\end{lemma}

\begin{lemma} Consider $k(r) \tau \bt dr$ where $\tau$ satisfies
(\ref{decloc20}) and (\ref{decloc30}).  Then
\begin{gather} \label{crosslap}
\begin{split}
\N^* \N (k \tau \bt dr) =&\ - 2 r^{-2} k \cn{\N}^{\sym} \tau\\
&\ + \left\{ \left[ -k'' - 3 r^{-1} k' + r^{-2} k \left( \cn{\N}^* \cn{\N} + 6
\right) \right] \tau \right\} \bt dr\\
&\ + 4 r^{-2} k \cn{\gd} \tau dr \otimes dr
\end{split}
\end{gather}
where
\begin{gather}
\cn{\N}^{\sym} \tau = \sum_i \til{\N}_{e_i} \tau \bt e_i^*
\end{gather}
\end{lemma}

\begin{lemma} Consider $l(r) \phi(x) dr \otimes dr$.  Then
\begin{gather} \label{vertlap}
\begin{split}
\N^* \N (l \phi dr \otimes dr) =&\ - 2 r^{-2} l \phi \cn{g} - 2 r^{-2} l
\cn{d} \phi \bt dr\\
&\ + (- l'' - 3 r^{-1} l' + r^{-2} l (\cn{\N}^* \cn{\N} + 6)) \gp dr \otimes
dr
\end{split}
\end{gather}
\end{lemma}

\begin{lemma} \label{divg} For tensors of the three types considered above, we
have the following formulas:
\begin{align}
\gd \left(f B \right) =&\ r^{-1} f \cn{\gd}(B) - r^{-1} f \cn{\tr}(B) dr\\
\gd \left( k \tau \bt dr \right) =&\ k' \tau + 4 r^{-1} k \tau + r^{-1} k
\left(\cn{\gd} \tau \right) dr\\
\gd \left( l \gp dr \otimes dr \right) =&\ \left(l' + 3 r^{-1} l \right) \gp
dr
\end{align}
\end{lemma}

\begin{prop} Consider $f(r) B$ where $B = \sum \eta_i \bt \eta_j$ and $\eta_i$
satisfies (\ref{decloc20}) and (\ref{decloc30}).  Then
\begin{gather} \label{horizlap2}
\begin{split}
\left(\N^* \N \right)^2 (fB) =&\ \left(f^{(4)} +  6 r^{-1} f^{(3)} - r^{-2}
f''\left(1 + 2 \cn{\N}^* \cn{\N} \right) \right) B\\
&\ - \left( r^{-3} f' \left(2 \cn{\N}^* \cn{\N} + 7 \right) \right) B + r^{-4}
f \left( \left(\cn{\N}^* \cn{\N} + 2 \right) \right)^2 B\\
&\ - 4 r^{-4} f \cn{\N}^{\sym} \cn{\gd}(B) + 4 r^{-4} f \cn{\tr}(B) \cn{g}\\
&\ + 4 \left\{ \left[ - r^{-2} f'' - r^{-3} f' + r^{-4} f \left( \cn{\N}^*
\cn{\N} + 2 \right) \right] \cn{\gd}(B) \right\} \bt dr\\
&\ + 8 r^{-4} f \cn{d} \cn{\tr}(B) \bt dr\\
&\ + 4 \left( r^{-2} f'' + r^{-3} f' - r^{-4} f \left( \cn{\N}^* \cn{\N} + 4
\right) \right) \cn{\tr}(B) dr \otimes dr\\
&\ + 8 r^{-4} f \cn{\gd} \cn{\gd}(B) dr \otimes dr
\end{split}
\end{gather}
\begin{proof} This will be a straightforward calculation using the three above
lemmas.  First of all using (\ref{horizlap}) we see that
\begin{gather} \label{hlsloc10}
\begin{split}
- \N^* \N \left( f'' B \right) =&\ \left( f^{(4)} + 3 r^{-1} f^{(3)} - r^{-2}
f'' \left(\cn{\N}^* \cn{\N} + 2 \right) \right) B\\
&\ - 2 r^{-2} f'' \cn{\gd}(B) \bt dr + 2 r^{-2} f'' \cn{\tr}(B) dr \otimes dr
\end{split}
\end{gather}
Next
\begin{gather} \label{hlsloc20}
\begin{split}
\N^* \N \left( -3 r^{-1} f' B \right) =&\ 3 \left( r^{-1} f^{(3)} + r^{-2} f''
- r^{-3} f' \left(\cn{\N}^* \cn{\N}+ 3 \right) \right) B\\
&\ - 6 r^{-3} f' \cn{\gd}(B) \bt dr + 6 r^{-3} f' \cn{\tr}(B) dr \otimes dr
\end{split}
\end{gather}
Now it is clear that $\til{\N}^* \til{\N} B$ satisfies (\ref{decloc20}) and
(\ref{decloc30}).  Thus apply (\ref{horizlap}) again to compute
\begin{gather}
\begin{split}
\N^* \N& \left( r^{-2} f \left(\cn{\N}^* \cn{\N} + 2 \right) B \right)\\
=&\ \left( - r^{-2} f'' + r^{-3} f' + r^{-4} f \left(\cn{\N}^* \cn{\N}+2
\right) \right) \left( \cn{\N}^* \cn{\N} + 2 \right)(B)\\
&\ + 2 r^{-4} f \left( \cn{\N}^* \cn{\N} - 2 \right)( \cn{\gd} B) \bt dr + 4
r^{-4} f \cn{d} \cn{\tr} B \bt dr\\
&\ - 2 r^{-4} f \left( \cn{\N}^* \cn{\N} + 2 \right) \cn{\tr}(B) dr \otimes dr
\end{split}
\end{gather}
Note that we have applied lemma \ref{S3calc} to commute $\cn{\gd}$ past
$\cn{\N}^* \cn{\N}$. Now, since $\eta$ satisfies (\ref{decloc20}) and
(\ref{decloc30}), so does $\cn{\gd} B$, thus we may apply (\ref{crosslap}) to
conclude
\begin{gather}
\begin{split}
\N^* \N &\left( 2 r^{-2} f \cn{\gd} (B) \bt dr \right)\\
=&\ -4 r^{-4} f \cn{\N}^{\sym} \cn{\gd}(B)\\
&\ + \left\{ \left[ - 2 r^{-2} f'' + 2 r^{-3} f' + 2 r^{-4} f \left( \cn{\N}^*
\cn{\N} + 6 \right) \right] \cn{\gd}(B) \right\} \bt dr\\
&\ + 8 r^{-4} f \cn{\gd} \cn{\gd}(B) dr \otimes dr
\end{split}
\end{gather}
Similarly we apply (\ref{vertlap}) to conclude
\begin{gather} \label{hlsloc30}
\begin{split}
\N^* \N &\left( -2 r^{-2} f \cn{\tr}(B) dr \otimes dr \right)\\
=&\ 4 r^{-4} f \cn{\tr}(B) \cn{g} + 4 r^{-4} f \cn{d} \cn{\tr}(B) \bt dr\\
&\ + \left( 2 r^{-2} f'' - 2 r^{-3} f' - 2 r^{-4} f \left( \cn{\N}^* \cn{\N}
+ 6 \right) \right) \cn{\tr}(B) dr \otimes dr
\end{split}
\end{gather}
Summing together (\ref{hlsloc10}) - (\ref{hlsloc30}) gives the result.
\end{proof}
\end{prop}

\begin{prop} Consider $k(r) \tau \bt dr$ where $\tau$ satisfies
(\ref{decloc20}) and (\ref{decloc30}).  Then
\begin{gather} \label{crosslap2}
\begin{split}
\left(\N^* \N \right)^2 &\left( k \tau \bt dr \right)\\
=&\ 4 \left( r^{-2} k'' + r^{-3} k' - r^{-4} k \left(\cn{\N}^* \cn{\N}+5
\right) \right) \left(\cn{\N}^{\sym} \tau \right)\\
&\ - 8 r^{-4} k \cn{\gd} \tau \cn{g}\\
&\ + \left( k^{(4)} + 6 r^{-1} k^{(3)} - r^{-2} k'' \left( 2 \cn{\N}^* \cn{\N}
+ 9 \right) - r^{-3} k' \left(2 \cn{\N}^* \cn{\N} + 15 \right) \right.\\
&\ \left. + r^{-4} k \left( \left( \cn{\N}^* \cn{\N} \right)^2 + 16 \cn{\N}^*
\cn{\N} + 28 \right) \right) \tau \bt dr\\
&\ - 12 r^{-4} k \cn{d} \cn{\gd} \tau \bt dr\\
&\ + 8 \left( - r^{-2} k'' - r^{-3} k' + r^{-4} k \left( \cn{\N}^* \cn{\N} + 6
\right) \right) \cn{\gd} \tau dr \otimes dr
\end{split}
\end{gather}
\begin{proof} We must start from the expression in (\ref{crosslap}) and apply
(\ref{horizlap}), (\ref{crosslap}), and (\ref{vertlap}) to the individual
terms.  First of all it is clear that if $\tau$ satisfies (\ref{decloc20}) and
(\ref{decloc30}) then so does $\cn{\N}^{\sym} \tau$, thus we may apply
(\ref{horizlap}) to conclude
\begin{gather} \label{clsloc10}
\begin{split}
\left(\N^* \N \right)& \left( -2 r^{-2} k \cn{\N}^{\sym} \tau \right)\\
=&\ \left(2 r^{-2} k'' - 2 r^{-3} k' - 2 r^{-4} k \left(\cn{\N}^* \cn{\N}+2
\right) \right) \left(\cn{\N}^{\sym} \tau \right)\\
&\ + 4 r^{-4} k \left(\cn{\N}^* \cn{\N} - \cn{d} \cn{\gd} - 2\right)\tau \bt dr
+ 8 r^{-4} k \cn{\gd} \tau dr \otimes dr
\end{split}
\end{gather}
where we have used the equations
\begin{gather*}
\cn{\gd} \left( \cn{\N}^{\sym} \tau \right) = - \cn{\N}^* \cn{\N} \tau + \cn{d}
\cn{\gd} \tau + 2 \tau, \qquad \cn{\tr} \left( \cn{\N}^{\sym} \tau \right) = 2
\cn{\gd}
\tau
\end{gather*}
from lemma \ref{S3calc}.  Now we apply (\ref{crosslap}) to conclude
\begin{gather}
\begin{split}
\left(\N^* \N \right)& \left(- k'' \tau \bt dr \right)\\
=&\  2 r^{-2} k'' \cn{\N}^{\sym} \tau\\
&\ + \left\{ \left[ k^{(4)} + 3 r^{-1} k^{(3)} - r^{-2} k''\left( \cn{\N}^*
\cn{\N} + 6 \right) \right] \tau \right\} \bt dr\\
&\ - 4 r^{-2} k'' \cn{\gd} \tau dr \otimes dr
\end{split}
\end{gather}
And similarly
\begin{gather}
\begin{split}
\left(\N^* \N \right)& \left(- 3 r^{-1} k' \tau \bt dr \right)\\
=&\  6 r^{-3} k' \cn{\N}^{\sym} \tau\\
&\ + 3 \left\{ \left[ r^{-1} k^{(3)} + r^{-2} k'' - r^{-3} k' \left(\cn{\N}^*
\cn{\N} + 7 \right) \right] \tau \right\} \bt dr\\
&\ - 12 r^{-3} k' \cn{\gd} \tau dr \otimes dr
\end{split}
\end{gather}
And again
\begin{gather}
\begin{split}
\left(\N^* \N \right)& \left( r^{-2} k \left( \cn{\N}^* \cn{\N} + 6 \right)
\tau \right) \bt dr\\
=&\ -2 r^{-4} k \left( \cn{\N}^* \cn{\N} + 8 \right) \cn{\N}^{\sym} \tau\\
&\ + \left\{ \left[ - r^{-2} k'' + r^{-3} k' + r^{-4} k \left( \cn{\N}^*
\cn{\N} + 6 \right) \right] \left( \cn{\N}^* \cn{\N} + 6 \right) \tau \right\}
\bt dr\\
&\ + 4 r^{-4} k \left(\cn{\N}^* \cn{\N} + 4 \right) \cn{\gd} \tau dr \otimes
dr
\end{split}
\end{gather}
where we commuted derivatives using lemma \ref{S3calc}.  Finally using
(\ref{vertlap}) with $l = r^{-2} k$ and $\phi = \cn{\gd} \tau$ we compute
\begin{gather} \label{clsloc20}
\begin{split}
\left(\N^* \N \right)& \left( 4 r^{-2} k \cn{\gd} \tau \right) dr \otimes dr\\
=&\ -8 r^{-4} k \cn{\gd} \tau \cn{g} - 8 r^{-4} k \cn{d} \cn{\gd} \tau \bt dr\\
&\ + 4 \left[ - r^{-2} k'' + r^{-3} k' + r^{-4} k \left( \cn{\N}^* \cn{\N} +
6 \right) \right] \left(\cn{\gd} \tau\right) dr \otimes dr
\end{split}
\end{gather}
Collecting together (\ref{clsloc10}) - (\ref{clsloc20}) gives the result.
\end{proof}
\end{prop}

\begin{prop} Consider $l(r) \phi(x) dr \otimes dr$.  Then
\begin{gather} \label{vertlap2}
\begin{split}
\left(\N^* \N \right)^2& \left( l \phi dr \otimes dr \right)\\
=&\ 4 \left( r^{-2} l'' + r^{-3} l' - r^{-4} l \left( \cn{\N}^* \cn{\N} + 4
\right) \right) \left( \phi \cn{g} \right)\\
&\ + 4 r^{-4} l \cn{\N}^{\sym} \cn{d} \phi\\
&\ + 4 \left( r^{-2} l'' + r^{-3} l' - r^{-4} l \left( \cn{\N}^* \cn{\N} + 8
\right) \right) \cn{d} \phi \bt dr\\
&\ + \left( l^{(4)} + 6 r^{-1} l^{(3)} - r^{-2} l''\left(2 \cn{\N}^* \cn{\N}
+ 9 \right) - r^{-3} l' \left( 2 \cn{\N}^* \cn{\N} + 15 \right) \right.\\
&\ \left. + r^{-4} l \left( \left(\cn{\N}^* \cn{\N} \right)^2 + 20 \cn{\N}^*
\cn{\N} + 48 \right) \right) \gp dr \otimes dr
\end{split}
\end{gather}

\begin{proof} We need to compute the action of the Laplacian of the right hand
side of (\ref{vertlap}).  First of all it is clear that $\phi \cn{g}$ satisfies
(\ref{decloc20}) and (\ref{decloc30}) thus we may apply (\ref{horizlap})
to conclude
\begin{gather} \label{vlsloc10}
\begin{split}
\left( \N^* \N \right)& \left( - 2 r^{-2} l \phi \cn{g} \right)\\
=&\ \left(2 r^{-2} l'' - 2 r^{-3} l' - 2 r^{-4} l \left( \cn{\N}^* \cn{\N}+2
\right) \right) \left(\phi \cn{g} \right)\\
&\ - 4 r^{-4} l \cn{d} \phi \bt dr + 12 r^{-4} l \phi dr \otimes dr\\
\end{split}
\end{gather}
where we used that $\cn{\gd} \cn{g} = 0$ and $\cn{\tr} \cn{g} = 3$.  Next we
note that $\cn{d} \phi$ satisfies (\ref{decloc20}) and (\ref{decloc30}) so
that we can apply (\ref{crosslap}) to conclude
\begin{gather}
\begin{split}
\left( \N^* \N \right)& \left( - 2 r^{-2} l \right) \cn{d} \phi \bt dr\\
=&\ 4 r^{-4} l \cn{\N}^{\sym} \cn{d} \phi\\
&\ + \left\{ \left[ 2 r^{-2} l'' - 2 r^{-3} l' - 2 r^{-4} l \left( \cn{\N}^*
\cn{\N} + 6 \right) \right] \cn{d} \phi \right\} \bt dr\\
&\ + 8 r^{-4} l \cn{\N}^* \cn{\N} \phi dr \otimes dr
\end{split}
\end{gather}
Where we have used the equation \textbf{$\cn{\N}^* \cn{\N} = - \cn{\gd} \cn{d}$}
Next using (\ref{vertlap}) we compute
\begin{gather}
\begin{split}
\left(\N^* \N \right)& \left(- l'' \phi dr \otimes dr \right)\\
=&\ 2 r^{-2} l'' \phi \cn{g} + 2 r^{-2} l'' \cn{d} \phi \bt dr\\
&\ + \left( l^{(4)} + 3 r^{-1} l^{(3)} - r^{-2} l'' \left( \cn{\N}^* \cn{\N}
+ 6 \right) \right) \phi dr \otimes dr
\end{split}
\end{gather}
And similarly
\begin{gather}
\begin{split}
\left( \N^* \N \right)& \left( - 3 r^{-1} l' \phi dr \otimes dr \right)\\
=&\ 6 r^{-3} l' \phi \cn{g} + 6 r^{-3} l' \cn{d} \phi \bt dr\\
&\ + 3 \left(r^{-1} l^{(3)} + r^{-2} l'' - r^{-3} l' \left( \cn{\N}^* \cn{\N}
+ 7) \right) \right) \phi dr \otimes dr
\end{split}
\end{gather}
Finally we apply (\ref{vertlap}) to conclude
\begin{gather} \label{vlsloc20}
\begin{split}
\left( \N^* \N \right)& \left(r^{-2} l \left( \cn{\N}^* \cn{\N} + 6 \right)
\phi dr \otimes dr \right)\\
=&\ - 2 r^{-4} l \left[\left( \cn{\N}^* \cn{\N} + 6 \right) \phi \right]
\cn{g} - 2 r^{-4} l \left[\left( \cn{\N}^* \cn{\N} + 8 \right) \cn{d} \phi
\right] \bt dr\\
&\ + \left( - r^{-2} l'' + r^{-3} l' + r^{-4} l \left( \cn{\N}^* \cn{\N} + 6
\right) \right) \left[\left( \cn{\N}^* \cn{\N} + 6 \right) \phi \right] dr
\otimes dr
\end{split}
\end{gather}
Collecting together (\ref{vlsloc10}) - (\ref{vlsloc20}) gives the result.
\end{proof}
\end{prop}

\section{The Linearized Equation}
In this section we derive the equation for the linearization of the Bach-flat
condition at a flat metric.  The equation is of course not strictly elliptic due
to the diffeomorphism invariance of the Bach-flat equation.  After restricting
to the case where the variation of the metric is trace-free and divergence-free,
we are reduced to the Biharmonic equation.  We can then use the decomposition of
the Bilaplacian computed in the previous section to classify solutions to this
equation and moreover compute the exact decay rates of solutions, which is the
key ingredient of our removal of singularities result.

\begin{prop} Given $g$ a flat metric, the linearization
of the Bach-flat equation at $g$ is equivalent to the equation
\begin{gather} \label{sloc5}
\begin{split}
0 =&\ \gD^2 h - 2 \gD \gd^* \gd h\\
&\ + \frac{1}{3} \left( \N^2 \gD \tr_g h - \left(\gD^2 \tr_g h \right) g
 + \left(\gD \gd^2 h \right) g + 2 \N^2 \gd^2 h \right)
\end{split}
\end{gather}
\begin{proof} Let $g(s)$ be a family of metrics such that $g(0) = g$ and
$\frac{\del}{\del s} g(s)_{|s = 0} = h$.  First of all since the metric is flat
it is clear that
\begin{gather*}
\left( \frac{1}{2} R^{kl} W_{ikjl} \right)' \cdot h = 0
\end{gather*}
As a consequence of the Bianchi identity, in dimension $4$ we have the
equation
\begin{gather*}
\N^k \N^l W_{ikjl} = \gD A_{ij} - \N^k \N_i A_{jk}
\end{gather*}
where $A_{ij}$ are the components of the Weyl-Schouten tensor
\begin{gather*}
A = \frac{1}{n - 2} \left( \Rc - \frac{1}{2(n - 1)} R \cdot g \right)
\end{gather*}
A further application of the Bianchi identity then yields
\begin{gather*}
\N^k \N^l W_{ikjl} = \frac{1}{2} \left( \gD \Rc_{ij} - \frac{1}{6} \gD R
g_{ij} - \frac{1}{3} \N_i \N_j R \right)
\end{gather*}
First we compute the linearization
\begin{gather*}
\begin{split}
\left( \gD \Rc \right)' \cdot h =&\ - \frac{1}{2} \gD \left( \gD h + \N^2
\tr h - 2 \gd^* \gd h \right)\\
\end{split}
\end{gather*}
and next
\begin{gather*}
\left( - \frac{1}{6} \gD R g \right)' \cdot h = \frac{1}{6} \gD \left(
\gD \tr h - \gd \gd h) \right) g
\end{gather*}
and finally
\begin{gather*}
\left( - \frac{1}{3} \N^2 R \right)' \cdot h = \frac{1}{3} \N^2
\left( \gD \tr h - \gd \gd h) \right)
\end{gather*}
The result follows.
\end{proof}
\end{prop}
Now, by a diffeomorphism gauge-fixing procedure which we will use in the
applications, we may assume that $h$ comes to us almost divergence-free.  In
fact we cannot prescribe that $h$ is exactly divergence free.  We will address
this difficulty later in this section.
Also, because of the conformal invariance of the Bach-flat equation, the trace
of the right hand side of (\ref{sloc5})
vanishes.  Thus there is no a-priori separate equation that the trace of $h$
must satisfy.  However, using this conformal invariance in the applications, we
can assume that our Bach-flat metric has constant scalar curvature, which in the
asymptotically flat regime corresponds to the equation
\begin{gather} \label{sloc7}
\gD \tr h = 0
\end{gather}
As it
turns out we can now remove the trace entirely by adding a term $L_X g$ for some
vector field $X$.  This corresponds to solving for a vector field $X$ satisfying
\begin{gather}
\gd L_X g = 0
\end{gather}
to preserve the divergence-free condition and
\begin{gather}
\tr \left( L_X g \right) = 2 \gd X = u
\end{gather}
where $u$ is a prescribed harmonic function.  It is shown in \cite{Tian} lines
4.11-4.20 that on a Ricci-flat cone one can indeed solve for this $X$.  Thus a
divergence-free solution to (\ref{sloc5}) satisfying (\ref{sloc7}) can be
written
as $L_X g + h$ where
\begin{gather} \label{sloc30}
\begin{split}
0 =&\ \gD^2 h\\
0 =&\ \tr_g h\\
0 =&\ \gd h
\end{split}
\end{gather}
We now proceed to analyse solutions to (\ref{sloc30}).  Let
$h = f(r) B + k(r) \tau \bt dr + l(r) \phi(x) dr \otimes dr$.  First
of all the equation $\tr h = 0$ reduces to the equation
\begin{gather} \label{sloc50}
f \cn{\tr}(B) + l \gp = 0
\end{gather}
Next, using lemma \ref{divg} and (\ref{sloc50}) the equation $\gd h = 0$
reduces to the equations
\begin{align} \label{sloc60}
r^{-1} f \cn{\gd}(B) + \left(k' + 4 r^{-1} k \right) \tau =&\ 0\\
\label{sloc70}
\left(l' + 4 r^{-1} l \right) \gp + r^{-1} k \left(\cn{\gd} \tau \right) =&\ 0
\end{align}
Now we would like to write down and simplify the vertical component of
equation (\ref{sloc30}).  Using propositions \ref{horizlap2}, \ref{crosslap2}
and \ref{vertlap2} we see that this equation gives
\begin{gather} \label{sloc95}
\begin{split}
0 =&\ \left( l^{(4)} + 6 r^{-1} l^{(3)} - r^{-2} l''\left(2 \cn{\N}^* \cn{\N}
+ 9 \right) - r^{-3} l' \left( 2 \cn{\N}^* \cn{\N} + 15 \right) \right.\\
&\ \left. + r^{-4} l \left( \left(\cn{\N}^* \cn{\N} \right)^2 + 20 \cn{\N}^*
\cn{\N} + 48 \right) \right) \gp\\
&\ + 8 \left( - r^{-2} k'' - r^{-3} k' + r^{-4} k \left( \cn{\N}^* \cn{\N} +
6 \right) \right) \cn{\gd} \tau\\
&\ + 4 \left( r^{-2} f'' + r^{-3} f' - r^{-4} f \left( \cn{\N}^* \cn{\N} + 4
\right) \right) \cn{\tr}(B)\\
&\ + 8 r^{-4} f \cn{\gd} \cn{\gd}(B)
\end{split}
\end{gather}
Now we want to use equations (\ref{sloc50}), (\ref{sloc60}) and (\ref{sloc70})
to simplify this further.  First of all using (\ref{sloc50}) it is clear that
\begin{gather} \label{sloc100}
\begin{split}
4 & \left(r^{-2} f'' + r^{-3} f' - r^{-4} f\left(\cn{\N}^* \cn{\N} + 4 \right)
\right) \cn{\tr}(B)\\
&\ \quad = - 4 \left( r^{-2} l'' + r^{-3} l' - r^{-4} l \left(\cn{\N}^*
\cn{\N} + 4 \right) \right) \gp
\end{split}
\end{gather}
Next using (\ref{sloc70}) we compute
\begin{gather} \label{sloc110}
\begin{split}
8 & \left( - r^{-2} k'' - r^{-3} k' + r^{-4} k\left( \cn{\N}^* \cn{\N} + 6
\right) \right) \cn{\gd} \tau\\
&\ \quad = 8 \left( r^{-1} l^{(3)} + 7 r^{-2} l'' - r^{-3} l'\left( \cn{\N}^*
\cn{\N} + 1 \right) - 4 r^{-4} l \left( \cn{\N}^* \cn{\N} + 6 \right) \right)
\gp
\end{split}
\end{gather}
Finally using both (\ref{sloc60}) and (\ref{sloc70}) we compute
\begin{gather} \label{sloc120}
\begin{split}
8 r^{-4} f \cn{\gd} \cn{\gd}(B) =&\ - 8 \left( r^{-3} k' + 4 r^{-4} k \right)
\cn{\gd} \tau\\
=&\ 8 \left(r^{-2} l'' + 9 r^{-3} l' + 16 r^{-4} l \right) \gp
\end{split}
\end{gather}
Thus plugging (\ref{sloc100}) - (\ref{sloc120}) into (\ref{sloc95}) gives the
equation
\begin{gather}
\begin{split}
0 =&\ \left( l^{(4)} + 14 r^{-1} l^{(3)} + r^{-2} l''\left(- 2 \cn{\N}^*
\cn{\N} + 51 \right) \right.\\
&\ \left. + r^{-3} l' \left( - 10 \cn{\N}^* \cn{\N} + 45 \right) + r^{-4} l
\left( \left( \cn{\N}^* \cn{\N} \right)^2 - 8 \cn{\N}^* \cn{\N} \right)
\right) \gp
\end{split}
\end{gather}
We would now like to do the same for the cross component of (\ref{sloc30}).
First of all from propositions \ref{horizlap2}, \ref{crosslap2} and
\ref{vertlap2} we see that this is the equation
\begin{gather} \label{sloc150}
\begin{split}
0 =&\ 4 \left[ - r^{-2} f'' - r^{-3} f' + r^{-4} f \left( \cn{\N}^* \cn{\N} +
2 \right) \right] \cn{\gd}(B) + 8 r^{-4} f \cn{d} \cn{\tr}(B)\\
&\ + \left( k^{(4)} + 6 r^{-1} k^{(3)} - r^{-2} k'' \left( 2 \cn{\N}^* \cn{\N}
+ 9 \right) - r^{-3} k' \left(2 \cn{\N}^* \cn{\N} + 15 \right) \right.\\
&\ \left. + r^{-4} k \left( \left( \cn{\N}^* \cn{\N} \right)^2 + 16 \cn{\N}^*
\cn{\N} + 28 \right) \right) \tau - 12 r^{-4} k \cn{d} \cn{\gd} \tau\\
&\ + 4 \left( r^{-2} l'' + r^{-3} l' - r^{-4} l \left( \cn{\N}^* \cn{\N} + 8
\right) \right) \cn{d} \phi
\end{split}
\end{gather}
Now, using (\ref{sloc60}) we see that
\begin{gather} \label{sloc160}
\begin{split}
4 &\left[ - r^{-2} f'' - r^{-3} f' + r^{-4} f \left( \cn{\N}^* \cn{\N} + 2
\right) \right] \cn{\gd}(B)\\
&\ \quad = 4 \left(r^{-1} k^{(3)} + 7 r^{-2} k'' + r^{-3} k'\left( - \cn{\N}^*
\cn{\N} + 3 \right) \right.\\
&\ \left. \quad\ - 4 r^{-4} k \left(\cn{\N}^* \cn{\N} + 2 \right) \right) \tau
\end{split}
\end{gather}
Next using (\ref{sloc50}) we make the simplification
\begin{gather}
8 r^{-4} f \cn{d} \cn{\tr}(B) = -8 r^{-4} l \cn{d} \gp
\end{gather}
And using (\ref{sloc70}) we simplify
\begin{gather} \label{sloc170}
- 12 r^{-4} k \cn{d} \cn{\gd} \tau = 12 \left( r^{-3} l' + 4 r^{-4} l \right)
\cn{d} \gp
\end{gather}
Plugging (\ref{sloc160}) - (\ref{sloc170}) into (\ref{sloc150}) gives
\begin{gather}
\begin{split}
0 =&\ \left( k^{(4)} + 10 r^{-1} k^{(3)} + r^{-2} k''\left(- 2 \cn{\N}^*
\cn{\N} + 19 \right) \right.\\
&\ \left. - k' \left( 6 \cn{\N}^* \cn{\N} + 3 \right) + r^{-4} k \left(
\left( \cn{\N}^* \cn{\N} \right)^2 - 4 \right) \right)
\tau\\
&\ + 4 \left( r^{-2} l'' + 4 r^{-3} l' + r^{-4} l\left( - \cn{\N}^* \cn{\N}
+ 2 \right) \right) \cn{d} \gp
\end{split}
\end{gather}
Finally we need to simplify the horizontal component of (\ref{sloc30}).  Again
using propositions \ref{horizlap2}, \ref{crosslap2} and \ref{vertlap2} we get
the equation
\begin{gather} \label{sloc200}
\begin{split}
0 =&\ \left(f^{(4)} +  6 r^{-1} f^{(3)} - r^{-2} f''\left(1 + 2 \cn{\N}^*
\cn{\N} \right) \right.\\
&\ \left. - r^{-3} f' \left(2 \cn{\N}^* \cn{\N} + 7 \right) + r^{-4} f
\left(\cn{\N}^* \cn{\N} + 2 \right)^2 \right) B\\
&\ - 4 r^{-4} f \cn{\N}^{\sym} \cn{\gd}(B) + 4 r^{-4} f \cn{\tr}(B) \cn{g}\\
&\ + 4 \left( r^{-2} k'' + r^{-3} k' - r^{-4} k \left(\cn{\N}^* \cn{\N}+5
\right) \right) \left(\cn{\N}^{\sym} \tau \right) - 8 r^{-4} k \cn{\gd} \tau
\cn{g}\\
&\ + 4 \left( r^{-2} l'' + r^{-3} l' - r^{-4} l \left( \cn{\N}^* \cn{\N} + 4
\right) \right) \left( \phi \cn{g} \right)\\
&\ + 4 r^{-4} l \cn{\N}^{\sym} \cn{d} \phi
\end{split}
\end{gather}
First using (\ref{sloc60}) we simplify
\begin{gather} \label{sloc210}
- 4 r^{-4} f \cn{\N}^{\sym} \cn{\gd}(B) = 4 \left( r^{-3} k' + 4 r^{-4} k
\right) \cn{\N}^{\sym} \tau
\end{gather}
Using (\ref{sloc50}) we simplify
\begin{gather}
4 r^{-4} f \cn{\tr}(B) \cn{g} = - 4 r^{-4} l \gp \cn{g}
\end{gather}
Next from (\ref{sloc70}) we simplify
\begin{gather} \label{sloc220}
-8 r^{-4} k \cn{\gd} \tau \cn{g} = 8 \left( r^{-3} l' + 4 r^{-4} l \right)
\gp \cn{g}
\end{gather}
Plugging (\ref{sloc210}) - (\ref{sloc220}) into (\ref{sloc200}) gives
\begin{gather}
\begin{split}
0 =&\ \left(f^{(4)} +  6 r^{-1} f^{(3)} - r^{-2} f''\left(1 + 2 \cn{\N}^*
\cn{\N} \right) \right.\\
&\ \left. - r^{-3} f' \left(2 \cn{\N}^* \cn{\N} + 7 \right) + r^{-4} f
\left(\cn{\N}^* \cn{\N} + 2 \right)^2 \right) B\\
&\ + 4 \left( r^{-2} k'' + 2 r^{-3} k' - r^{-4} k \left(\cn{\N}^* \cn{\N} +1
\right) \right)
\left( \cn{\N}^{\sym} \tau \right)\\
&\ + 4 \left( r^{-2} l'' + 3 r^{-3} l' + r^{-4} l \left( - \cn{\N}^* \cn{\N} +
3 \right) \right) \gp \cn{g} + 4 r^{-4} l \cn{\N}^{\sym} \cn{d} \gp
\end{split}
\end{gather}

\begin{lemma} \label{syslemma} The system of equations (\ref{sloc30}) is
equivalent to
\begin{gather} \label{sysloc1}
\begin{split}
0 =&\ \left(f^{(4)} +  6 r^{-1} f^{(3)} - r^{-2} f''\left(1 + 2 \cn{\N}^*
\cn{\N} \right) \right.\\
&\ \left. - r^{-3} f' \left(2 \cn{\N}^* \cn{\N} + 7 \right) + r^{-4} f
\left(\cn{\N}^* \cn{\N} + 2 \right)^2 \right) B\\
&\ + 4 \left( r^{-2} k'' + 2 r^{-3} k' - r^{-4} k \left(\cn{\N}^* \cn{\N} +1
\right) \right)
\left( \cn{\N}^{\sym} \tau \right)\\
&\ + 4 \left( r^{-2} l'' + 3 r^{-3} l' + r^{-4} l \left( - \cn{\N}^* \cn{\N} +
3 \right) \right) \gp \cn{g} + 4 r^{-4} l \cn{\N}^{\sym} \cn{d} \gp
\end{split}
\end{gather}
\begin{gather} \label{sysloc2}
\begin{split}
0 =&\ \left( k^{(4)} + 10 r^{-1} k^{(3)} + r^{-2} k''\left(- 2 \cn{\N}^*
\cn{\N} + 19 \right) \right.\\
&\ \left. - r^{-3} k' \left( 6 \cn{\N}^* \cn{\N} + 3 \right) + r^{-4} k \left(
\left( \cn{\N}^* \cn{\N} \right)^2 - 4 \right) \right)
\tau\\
&\ + 4 \left( r^{-2} l'' + 4 r^{-3} l' + r^{-4} l\left( - \cn{\N}^* \cn{\N}
+ 2 \right) \right) \cn{d} \gp
\end{split}
\end{gather}
\begin{gather} \label{sysloc3}
\begin{split}
0 =&\ \left( l^{(4)} + 14 r^{-1} l^{(3)} + r^{-2} l''\left(- 2 \cn{\N}^*
\cn{\N} + 51 \right) \right.\\
&\ \left. + r^{-3} l' \left( - 10 \cn{\N}^* \cn{\N} + 45 \right) + r^{-4} l
\left( \left( \cn{\N}^* \cn{\N} \right)^2 - 8 \cn{\N}^* \cn{\N} \right)
\right) \gp
\end{split}
\end{gather}
\begin{align} \label{sysloc4}
0 =&\ r^{-1} f \cn{\gd}(B) + \left(k' + 4 r^{-1} k \right) \tau\\
\label{sysloc5}
0 =&\ \left(l' + 4 r^{-1} l \right) \gp + r^{-1} k \left(\cn{\gd} \tau \right)
\end{align}
\end{lemma}

We now analyse all solutions to this system of equations.  We notice that for
any solution, $l$ satisfies the determined ODE in (\ref{sysloc3}).  Thus we can
compute the solutions to (\ref{sysloc3}) in terms of the eigenvalues of the
Laplacian on
$S^3$ acting on functions.  Once these are determined we notice that the
equations (\ref{sysloc4}) and (\ref{sysloc5}) show that both $f$ and $k$ will
have the same decay rate as $l$.  Once this is done we restrict to the case
where $l = 0$.  Here we notice that $k$ satisfies the now determined ODE in
(\ref{sysloc2}).  We can classify these solutions in terms of the eigenvalues of
the Laplacian on $S^3$ acting on one-forms.  Once again we can use equation
(\ref{sysloc4}) to conclude that $f$ has the same decay rate as $k$.  Finally we
restrict to the case where $l = k = 0$.  Then it is clear that $f$ satisfies the
determined ODE in (\ref{sysloc1}), and classify solutions in terms of the
eigenvalues of the Laplacian on $S^3$ acting on traceless symmetric two-tensors.
 We now make this rigorous in a series of lemmas.

\begin{lemma} \label{ODESolv1} The solutions to (\ref{sloc30}) satisfying $l = k
= 0$ are
\begin{gather*}
h = \sum_i r^{b_i} B_i
\end{gather*}
where either $b_i > 0$ or
\begin{gather*}
b_i \leq -2
\end{gather*}
\begin{proof}
Since $k = l = 0$ we have $\cn{\tr} B = \cn{\gd} B = 0$.  Say $\gl$ is an
eigenvalue for the Laplacian of
$S^3$ acting on traceless symmetric two-tensors, and suppose $B_\gl$ is in the
eigenspace of $\gl$.  Equation (\ref{sysloc1}) then reduces to an ODE with
solutions
\begin{gather}
r^{\pm 1 \pm \sqrt{\gl + 3}} B_\gl.
\end{gather}
From lemma \ref{evalcalc} we know that the smallest eigenvalue is $6$, so the
result follows.
\end{proof}
\end{lemma}

\begin{lemma} \label{ODESolv2} The solutions to (\ref{sloc30}) satisfying $l
= 0$ are
\begin{gather*}
h = L_X g_0 + \sum_i r^{b_i} \left(B_i + \tau_i \bt dr \right)
\end{gather*}
where either $b_i > 0$ or
\begin{gather*}
b_i \leq -2
\end{gather*}
\begin{proof} Note that since $l = 0$ equation (\ref{sysloc5}) implies in that
$\cn{\gd} \tau = 0$.  The eigenvalues for the Laplacian acting on $1$-forms on
$S^3$ are given by \cite{evals}
\begin{gather}
a_j := (j + 1)(j + 3)
\end{gather}
Say $\tau_j$ is in the eigenspace of $a_j$.  Then (\ref{sysloc2}) reduces to an
ODE with solutions
\begin{gather}
r^{ -2 \pm \sqrt{a_j + 2}} \tau_j
\end{gather}
We point out that the decay rates are all less than $-2$.  Also, since the decay
rate $-4$ does not occur, the expression $k' + 4 r^{-1} k$ never vanishes, thus
by equation (\ref{sysloc4}) we see that $f = c r^{-2 \pm \sqrt{a_j + 2}}$, and
the result follows.
\end{proof}
\end{lemma}

\begin{lemma} \label{ODESolv3} The solutions to (\ref{sloc30}) satisfying $l
\neq 0$ are
\begin{gather*}
h = L_X g_0 + \sum_i r^{b_i} \left(B_i + \tau_i \bt dr + \gp dr \otimes dr
\right)
\end{gather*}
where either $b_i > 0$ or
\begin{gather*}
b_i \leq -2
\end{gather*}
\begin{proof} The eigenvalues of the Laplacian acting on functions on $S^3$ are
\begin{gather} \label{adef}
a_j := j(j + 2)
\end{gather}
Suppose $\gp_j$ is in the eigenspace of $a_j$.  Equation (\ref{sysloc3}) reduces
to an ODE with solutions
\begin{gather}
l(r) = r^{-2 \pm 1 \pm \sqrt{a_j + 1}}
\end{gather}
In the case where $l(r) = r^{-3 \pm \sqrt{a_j + 1}}$, the solution to the whole
system is given by a Lie derivative term, specifically $L_X g_0$ where
\begin{gather*}
X_* = r^{-1 \pm \sqrt{a_j + 1}} d \gp_{j} + \left(-1 \pm \sqrt{a_j + 1} \right)
r^{-2 \pm \sqrt{a_j + 1}} \gp_j dr.
\end{gather*}

For the solutions $l(r) = r^{-1 \pm \sqrt{a_j + 1}}$ we note that for $j > 0$
the decay rates are always less than or equal to $-2$.  In the case $j = 0$,
i.e. where $\gp = c$, consider the radially parallel
solution, i.e. where $l(r) = c$.  Using equation (\ref{sysloc5}) we see that $k$
is also a nonzero constant.  Using this equation
(\ref{sysloc3}) reduces to $\left( \left(\cn{\N}^* \cn{\N} \right)^2 - 4 \right)
\tau = 0$ so we see that this solution does not in fact occur since
$2$ is not an eigenvalue of the Laplacian on $S^3$ acting on one-forms.

In any of the cases above, using equation (\ref{sloc50}) we see that $f$ will
have the same decay rate as $l$.  Given this, equation (\ref{sysloc4}) implies
that $k$ must also have the same decay rate.
\end{proof}
\end{lemma}

Together the above lemmas and earlier calculations prove the following
proposition.

\begin{prop} \label{linearsolns} On a flat cone, solutions of (\ref{sloc5})
satisfying $\gd h = 0$ and $\gD \tr h = 0$ can be written uniquely as a sum
\begin{gather}
h = L_X g_0 + \sum_i r^{b_i} \left(B_i + \tau_i \bt dr + \gp dr \otimes dr
\right)
\end{gather}
where either $b_i > 0$ or
\begin{gather} \label{growths}
b_i \leq -2
\end{gather}
\end{prop}

We hasten to point out that we will \emph{not} in fact be able to guarantee the
divergence-free condition.  This is due to the presence of certain eigenvalues
of the Laplace-Beltrami operator acting on one-forms.  Thus we follow the
technique used in \cite{Tian} and consider a modified divergence-free condition.
 In particular, fix $t \neq 0$ with $\brs{t}$ very small.  Let
\begin{gather} \label{divtdef}
\gd_t = \gd - t i_{r^{-1} \frac{\del}{\del r}}
\end{gather}
We will be able to prescribe that $\gd_t h = 0$ for $t$ arbitrarily small.
Given this, we define the following modified equation
\begin{gather} \label{modeqn}
\begin{split}
0 =&\ \gD^2 h - 2 t \gD \gd^* i_{r^{-1} \frac{\del}{\del r}} h\\
&\ + \frac{1}{3} \left(\N^2 \gD \tr_g h - \left(\gD^2 \tr_g h \right) g + t
\left(\gD \gd i_{r^{-1} \frac{\del}{\del r}} h \right) g + 2 t \N^2 \gd
i_{r^{-1} \frac{\del}{\del r}} h \right)\\
:=&\ \PP(h)
\end{split}
\end{gather}

Following the analysis of proposition \ref{linearsolns} we can write the
solutions satisfying $\gD \tr h = 0$ as growth and decay solutions where the
rates are perturbations of
those calculated above.  In particular we can write
\begin{gather}
h = r^{\gb_i} T_i
\end{gather}
where $T_i$ is some symmetric bilinear form and $\{T_i \}$ are orthonormal with
respect to the inner product
\begin{gather} \label{normdef1}
<<h_1, h_2>> = r^{-3} \int_{(r, S^3)} \left<h_1, h_2 \right> dV_{S^3}
\end{gather}
We also have the following corollary.

\begin{cor} For $t$ sufficiently small, there are no radially parallel solutions
of (\ref{modeqn}) satisfying $\gD \tr h = 0$.
\begin{proof} We note that the radially parallel solutions found above all had a
$dr$ component, and so in the perturbed equation are no longer radially
parallel.  Thus the only possibility would be $f(r) B$ where $f(r)$ is a
constant function and $B$ is trace and divergence-free.  However, no such
solution occurs according to our analysis above.  Thus the corollary follows.
\end{proof}
\end{cor}

So, according to the above results we can decompose any solution satisfying $\gD
\tr h = 0$ as
\begin{gather} \label{evaldecomp}
h = h_\up + h_\dn
\end{gather}
where $h_\up$ are the solutions with positive growth rate and likewise $h_\dn$
are solutions which decay in $r$.  Now let
\begin{gather}
\gb = \min_{i} \brs{\gb_i} > 0
\end{gather}
We can now state the first of our decay estimates.  We start with
some notation.  For a fixed annulus $A_{a,b}(p)$ we have the norm
\begin{gather} \label{normdef2}
|||h|||_{a,b} = \int_a^b ||h||^2 r^{-1} dr
\end{gather}
where the norm $||\ ||$ is defined in \ref{normdef1}
This norm is defined so that if $w = a^{-2} \psi_a^*(h)$, where $\psi_a$ is
the natural scaling map $\psi_a(x,r) = (x,ar)$, then
\begin{gather}
|||h|||_{a,La} = |||w|||_{1,L}
\end{gather}
\begin{cor} Given $0 < \gb' < \gb$ there exists $l$ such that for all $a > 0$
and $L \geq l$
\begin{gather}
|||h_\up|||_{La, L^2 a} \geq L^{\gb'} |||h_\up|||_{a,La}\\
|||h_\dn|||_{La, L^2 a} \leq L^{-\gb'}|||h_\dn|||_{a,La}
\end{gather}
\end{cor}

\begin{lemma} \label{lineargrowth} Given $0 < \gb' < \gb$ there exists $L$ such
that if $k$ is a solution to (\ref{modeqn}) satisfying $\gd_t k = 0$ and $\gD
\tr k = 0$ then if
\begin{gather}
|||k|||_{L,L^2} \geq L^{\gb'} |||k|||_{1,L}
\end{gather}
then
\begin{gather}
|||k|||_{L^2,L^3} \geq L^{\gb'} |||k|||_{L,L^2}
\end{gather}
and if
\begin{gather}
|||k|||_{L^2,L^3} \leq L^{-\gb'} |||k|||_{L,L^2}
\end{gather}
then
\begin{gather}
|||k|||_{L,L^2} \leq L^{-\gb'} |||k|||_{1,L}
\end{gather}
\end{lemma}

\section{Asymptotic Curvature Decay}
In this section we give the proof of theorem \ref{decaythm}.  The proof will
follow the techniques used in \cite{Tian}.  We first recall the
definition of certain norms introduced in \cite{Tian}.  Let $A_u$ denote the
natural action of the scaling $\psi_u$ on tensors of type $(p,q)$, i.e.
\begin{gather} \label{Adef}
A_u = \left(\psi_{u^{-1}} \right)_* \otimes \dots \otimes \left( \psi_{u^{-1}}
\right)_* \otimes \psi_u^* \otimes \dots \otimes \psi_u^*
\end{gather}
Given $T$ a tensor of type $(p,q)$ we have
\begin{gather*}
\nmm{T(u,x)}{k,\ga}{0} = \nmm{u^{p - q} A_u T(1,x)}{k,\ga}{0}
\end{gather*}
where the norm on the right is the $C^{k, \ga}$-norm with respect to $g_0$ at
the point $(1,x)$.  Using this we define
\begin{gather}
\nmm{T}{k, \ga}{0} = \sup_{(u,x)} \nmm{T(u,x)}{k,\ga}{0}
\end{gather}
And more generally
\begin{gather}
\nmm{T}{k,\ga}{l} = \nmm{r^{-l} T}{k,\ga}{0}
\end{gather}
Given $T$ a tensor of type $(p,q)$ we write $T \in \TTs{p,q}{k,\ga}{l}$ if
$\nmm{T}{k,\ga}{l} < \infty$.

We will make use of a gauge-fixing theorem (\cite{Tian} Theorem 3.1).  For a
given cone with cone-point $p$, let $A_{c,d} = \{ (r,x) | c < r < d \}$.
\begin{prop} \label{gaugefixthm} Fix $t \neq 0$.  There exists $\gk(t, k)$ such
that
if $(C(N^{n-1}), g_0)$ is a Ricci-flat cone and $g$ is a metric on
$A_{c,d}(p) \subset C(N^{n-1})$ where $d / c \geq 2$ such that
\begin{gather}
\brs{g - g_0}_{k,\ga;0} < \gk(t,k)
\end{gather}
then there exists a diffeomorphism $\gp: A_{c,d}(p) \to A_{c,d}(p)$ such that
\begin{gather}
\gp^*g \in \TTs{0,2}{k,\ga}{0}
\end{gather}
and
\begin{gather}
\gd_t \left( \gp^* g - g_0 \right) = 0
\end{gather}
\end{prop}
Since we have assumed our given Bach-flat metric is asymptotically flat,
we can use this theorem to construct gauges relative to the flat metric on the
cone in the asymptotic regime.  Specifically, let $g_0$ be the flat metric on
$\mathbb R^4$ and let $g$ be a given Bach-flat metric on this cone.  Consider
afixed annulus $A_{c,d}$ and let
\begin{gather}
h = \gp^* g - g_0
\end{gather}
where $\gp$ is from proposition \ref{gaugefixthm} with respect to this fixed
annulus.  Then we note that $h$
satisfies
\begin{gather}
\gd_t h = 0
\end{gather}
and further
\begin{gather}
B_{g_0 + h} - B_{g_0} = 0
\end{gather}
which is a nonlinear elliptic equation on $h$ which we can think of as a
perturbation of the linearized deformation equation for $h$ small.  The
following lemma makes this precise.
\begin{lemma} \label{approxlemma}
\begin{gather}
\begin{split}
B_{g_0 + h} - B_{g_0} =&\ \PP(h) + F(h, g_0)
\end{split}
\end{gather}
where
\begin{gather} \label{Fdef}
\nmm{F}{k,\ga}{0} \leq C \sum_{p + q + r + s = 4} \nmm{\N^p h}{k - p,\ga}{-p}
\nmm{\N^q h}{k - q,\ga}{-q} \nmm{\N^r h}{k - r,\ga}{-r}
\nmm{\N^s h}{k - s,\ga}{-s}
\end{gather}
\end{lemma}
Also, since the metric $g$ is scalar-flat, we get the equation
\begin{gather}
R_{g_0 + h} - R_{g_0} = 0
\end{gather}
from which we derive a similar lemma.
\begin{lemma} \label{approxlemma2}
\begin{gather}
R_{g_0 + h} - R_{g_0} = - \gD \tr h - t \gd \left(i_{r^{-1} \frac{\del}{\del r}}
h \right) + F'(h,g_0)
\end{gather}
where
\begin{gather} \label{F'def}
\nmm{F'}{k,\ga}{0} \leq C \sum_{p + q = 2} \nmm{\N^p h}{k - p,\ga}{-p} \nmm{\N^q
h}{k - q,\ga}{-q}
\end{gather}
\end{lemma}

\begin{lemma} \label{ellipticest} There exists a small constant
$\chi(\Lambda, K)$ such that if $\nmm{h}{k,\ga}{0} \leq \chi$ then
\begin{gather}
\nmm{h(2a)}{k,\ga}{0} \leq c(\Lambda, K, k) ||| h |||_{a, 4a}
\end{gather}
\begin{proof}  Using the scale-invariance of the respective norms, it suffices
to consider the case $a = 1$.  The result follows in this case by standard
elliptic theory.
\end{proof}
\end{lemma}

Now let $\pi$ denote orthogonal projection onto $\ker \left( \PP(h) \right)$
with respect to the inner product defining the norm $|||\ |||_{La, L^2 a}$.
Similar to (\ref{evaldecomp}) we can write
\begin{gather}
\pi h = \left(\pi h \right)_\up + \left( \pi h \right)_\dn
\end{gather}

\begin{prop} \label{growthprop} There exists $\chi = \chi(\Lambda, K) > 0$
and $L$ large so that if $\nmm{h}{k,\ga}{0} < \chi$ then if
\begin{gather} \label{growthhyp}
|||h|||_{La,L^2a} \geq L^{\gb'}|||h|||_{a,La}
\end{gather}
then
\begin{gather} \label{growthconc}
|||h|||_{L^2a,L^3a} \geq L^{\gb'}|||h|||_{La,L^2a}
\end{gather}
and if
\begin{gather} \label{decayhyp}
|||h|||_{L^2a,L^3a} \leq L^{-\gb'}|||h|||_{La,L^2a}
\end{gather}
then
\begin{gather} \label{decayconc}
|||h|||_{La,L^2a} \leq L^{-\gb'}|||h|||_{a,La}
\end{gather}
Moreover, at least one of (\ref{growthhyp}) or
(\ref{decayhyp}) must hold.
\begin{proof} Suppose one had a sequence of gauges $\phi_i$ and solutions
$h_i$ where $\nmm{h_i}{k,\ga}{0} \to 0$ but none of the assertions of the
proposition hold.  By rescaling our solutions and using a compactness
argument (see \cite{Tian} lemma 5.22) we can produce a solution to
(\ref{modeqn}) satisfying $\gD \tr h = 0$ which contradicts lemma
\ref{lineargrowth}.
\end{proof}
\end{prop}

\noindent{\textbf{Proof of Theorem \ref{decaythm}}} We will adopt the notation
used in this section.  In particular choose $\chi$ and $L$
as in proposition \ref{growthprop}.  Let $g_0$ denote the standard flat metric
on $\mathbb R^4$.  Since we already know that the metric $g$ is ALE of order
$0$,
there exists $\Psi : \mathbb R^4 \to X$ so that
\begin{gather} \label{mainproof10}
\nmm{\Psi^* g - g_0}{k,\ga}{0} = o(r)
\end{gather}
Now consider the sequence of annuli $A_{L^i a, L^{i+1} a}(0)$.  Using
(\ref{mainproof10}) we may apply proposition \ref{gaugefixthm} and choose a
$\gd_t$-free gauge for $\Psi^* g$ with respect to $g_0$, $\gp_i$, on each
annulus.  By pushing forward this is equivalent to finding a sequence of flat
metrics $g_i$ where $\Psi^* g - g_i$ is $\gd_t$-free with respect to $g_i$.
\begin{gather*}
h_i := \Psi^* g - g_i
\end{gather*}
We note that
\begin{gather*}
\lim_{i \to \infty} |||h_i|||_{L^i a,L^{i+1} a} = 0.
\end{gather*}
Indeed, if this were not the case then inductively applying proposition
\ref{growthprop}, and in particular using that (\ref{growthhyp}) implies
(\ref{growthprop}) we can contradict (\ref{mainproof10}).  Since for any of the
flat metrics $g_i$, there are no radially parallel solutions to the linearized
deformation equation, we can conclude that (\ref{decayconc}) holds for all $i$.

Also, by passing to a subsequence, we can assume that for some flat metric
$g_\infty$ we have
\begin{gather}
\lim_{j \to \infty} \brs{g_i - g_\infty}_{k,\ga';0} = 0 \quad \left(\ga' < \ga
\right)
\end{gather}
Using this together with (\ref{decayconc}) and lemma \ref{ellipticest} we
conclude that
\begin{gather} \label{basicdecay}
\nmm{\Psi^*g - g_\infty}{k,\ga}{0} \leq c r^{-\gb'}
\end{gather}
This proves that $g$ is ALE of order $\gb'$.

To prove the full claim we will use an inductive procedure choosing better and
better gauges.  First of all we point out that the decay estimate
(\ref{basicdecay}) in fact suffices to find a completely divergence-free global
gauge
$\gp$ (Remark 3.23 \cite{Tian}).  So, let $h = \phi^* g - g_0$ be
divergence-free.
First we get a better estimate on the trace of $h$.  In particular, by lemma
\ref{approxlemma2}
we have
\begin{gather}
\gD \tr h = F'(h)
\end{gather}
where the estimate (\ref{basicdecay}) implies
\begin{gather}
F'(h) \in \TTs{0,2}{k-2,\ga}{-2 \gb - 2}
\end{gather}
Applying the Greens function we can conclude the existence of a function $f$ so
that
\begin{gather}
\gD f = F'(h)
\end{gather}
where
\begin{gather}
f \in \TTs{0,2}{k,\ga}{-2\gb}
\end{gather}
In particular this implies the equation
\begin{gather}
\gD \left( \tr h - f \right) = 0
\end{gather}
As a harmonic function on $\mathbb R^4$, we know that the decay rate is at least
$-2$, thus
\begin{gather}
\tr h - f \in \TTs{0,2}{k,\ga}{-2}
\end{gather}
This argument already can be inducted to show
\begin{gather}
\tr h \in \TTs{0,2}{k,\ga}{2 - \gd}
\end{gather}
for any $\gd > 0$.  We now show how to get the rest of the estimate.
By lemma \ref{approxlemma} we have
\begin{gather}
\PP(h) = F(h)
\end{gather}
where $\PP$ now refers to the case $t = 0$.  Using (\ref{basicdecay}) we
conclude
\begin{gather}
F(h) \in \TTs{0,2}{k-4,\ga}{-4\gb - 4}
\end{gather}
Again using a Greens function, this decay is sufficient to conclude the
existence of $h_1$ such that
\begin{gather}
\PP h_1 = F(h)
\end{gather}
where
\begin{gather}
h_1 \in \TTs{0,2}{k,\ga}{-4 \gb}
\end{gather}
Thus $\PP(h - h_1) = 0$.  Given our estimate on the trace of $h$ we can assume
that $h - h_1$ is trace-free.  We note that $\PP$ acting on the traceless piece
of $h - h_1$ also vanishes.  Using proposition \ref{linearsolns} we can write
\begin{gather}
h - h_1 = L_X g_0 + h_2
\end{gather}
where
\begin{gather}
h_2 \in \TTs{0,2}{k,\ga}{- \min \{4 \gb,2\}}
\end{gather}
and also
\begin{gather}
X \in \TTs{1,0}{k+1,\ga}{1 - \gb}
\end{gather}
Let $K_X$ be the diffeomorphism generated by taking the flow of $X$ to time
$1$.  It is clear from the above estimates that
\begin{gather}
K_X^*(g_0) - L_X g_0 \in \TTs{0,2}{k,\ga}{- \min \{4 \gb, 2 \}}
\end{gather}
Furthermore putting together the above estimates we get
\begin{gather}
K_{-X}^* \phi^* g - g_0 \in \TTs{0,2}{k,\ga}{- \min \{4 \gb, 2 \}}
\end{gather}
If $4 \gb \geq 2$ we are done.  If not, start over with
$h = K_{-X}^* \phi^* g - g_0$ and proceed as above with $\gb$ replaced by
$4 \gb - \ge$ for very small $\ge > 0$.  It is clear that by induction the
result follows. \qed

\vspace{0.2in}
\noindent{\emph{Remark:}} The result holds assuming harmonic curvature instead
of the
Bach-flat condition.  In paticular, using the Bianchi identity equation
(\ref{harmoniccurvature}) implies an equation of the form $\gD \Rc = \Rm * \Rm$.
 The linearized deformation equation is the same as the one analyzed in section
3, and so one can apply the argument in theorem \ref{decaythm} to conclude that
the metric is in fact ALE of order $\tau$ for any $\tau < 2$.

\section{Proof of Theorem \ref{maintheorem}}
Let $\left(X, d, x \right)$ and $B(x, 1) \backslash \{x \}$ be as in the
statement of the theorem.
We already know by theorem 1.1 of \cite{TVBF} that $B(x, 1)$ has a
$C^0$-orbifold structure at $x$.  Thus in particular, lifting to the universal
cover we may suppose $x = 0$ and $B(x, 1) \subset D(0,1)$ where $D$ denotes the
distance ball in the Euclidean metric.  Let $g$ be this metric which is smooth
away from the origin and $C^0$ on $B(0,1)$.  Fix a constant $s < 1$, and
consider a sequence of annuli in the metric $g$
\begin{gather*}
A_i := A\left(s^{i+1}, s^{i} \right)
\end{gather*}
Now let
\begin{gather}
\begin{split}
\Phi&: \mathbb R^4 \backslash \{0 \} \to \mathbb R^4 \backslash \{0 \}\\
\Phi& (r, \gt) = \left( \frac{1}{r}, \gt \right)
\end{split}
\end{gather}
be the inversion through the unit sphere and let $A'_i = \Phi^* A_i$, and $g' =
\Phi^* g$.  Let $\rho_i$ denote the distance function of the metric $g'$ on the
annulus $A'_i$.  Since we have assumed a bound on the Sobolev constant, the
$L^2$-norm of curvature, and that the first Betti number is bounded, we may
apply the analysis of section 4 in \cite{TVBF}.  In particular we may conclude
that $\rho_i \to r$ as $i \to \infty$.  In particular the annuli $A'_i$ are
approaching the standard annuli in $\mathbb R^4$.  Thus the asymptotic analysis
of section 4 above applies so we conclude that this metric is ALE of order
$\tau$ for any $\tau < 2$.  Pulling this estimate back to $B$ using the
spherical inversion gives
\begin{gather} \label{mploc10}
g = \gd_{ij} + \OO(r^{2 - \gd}).
\end{gather}
For $\gd$ suficiently small this implies that the metric has a $C^{1,\ga}$
extension through the origin.  Using results from \cite{DeTurck} this implies
the existence of a harmonic coordinate around the origin.  A simple
computation using the Bianchi identity shows that the Bach-flat condition
implies that the curvature satisfies an equation of the form
$\gD \Rc = \Rm * \Rc$.  Thus we view the Bach-flat equation in harmonic
coordinates as the system
\begin{gather} \label{mploc30}
\gD \Rc = \Rm * \Rc\\ \label{mploc40}
\gD g =\Rc + Q(g, \del g)
\end{gather}
Using the curvature estimate $\brs{\Rm} = \OO(r^{-\gd})$ it follows from
equation (\ref{mploc30}) that $\Rc \in W^{2,p}$ for any $p$.  Using equation
(\ref{mploc40}) it is clear that $g \in W^{3,p}$.  This allows us to
bootstrap and conclude that $g \in C^{\infty}$.
\qed

\vspace{0.2in}
\noindent{\emph{Remark:}} Again we point out that the result holds assuming
harmonic
curvature instead of the
Bach-flat condition.  The crucial step in the above argument is the curvature
decay rate, and we mentioned after the proof of theorem \ref{decaythm} that
metrics with harmonic curvature satisfy this estimate.

\section{Appendix: Analysis on the Three Sphere}

\begin{lemma} \label{evalcalc} Let $(S^3, g)$ be the round three-sphere.  The
smallest
eigenvalue of the rough Laplacian acting on traceless divergence-free
symmetric two-tensors is $6$.
\begin{proof} Using the identification of $S^3$ with $SU(2)$ we consider the
standard global left-invariant Milnor frame $X^1, X^2, X^3$ with
structure constants
\begin{gather*}
C_{23}^1 = C_{31}^2 = C_{12}^3 = -2
\end{gather*}
If $e_1, e_2, e_3$ denotes the corresponding coframe, we have the equation
\begin{gather} \label{eigloc1}
\N_{X^i} e_j = - e_k \gs(ijk)
\end{gather}
where $\gs(ijk)$ denotes the sign of the permutation $(ijk)$ and is zero if
any of $i, j$ and $k$ are equal.  Using this global frame we may write any
traceless symmetric two-tensor as
\begin{gather*}
B_{ij} e_i \bt e_j
\end{gather*}
and compute
\begin{gather*}
\begin{split}
\N^* \N B =&\ - B_{kl,ii} e_k \bt e_l - 2 B_{kl,i} \left( \N_{X^i} e_k \bt e_l
 + e_k \bt \N_{X^i} e_l \right)\\
 &\ - B_{kl} \left( \N_{X_i} \N_{X_i} e_k \bt e_l + 2 \N_{X^i} e_k \bt
 \N_{X^i} e_l + e_k \bt \N_{X^i} \N_{X^i} e_l \right)
\end{split}
\end{gather*}
Now a basic calculation using (\ref{eigloc1}) shows that
\begin{gather*}
\N_{X^i} \N_{X^i} e_k = - 2 e_k
\end{gather*}
and similarly, using that $B$ is traceless it is easy to compute
\begin{gather*}
2 B_{kl} \N_{X^i} e_k \bt \N_{X^i} e_l = - 2 B_{kl} e_k \bt e_l
\end{gather*}
Using these calculations and the fact that the smallest eigenvalue must
occur when the coefficient functions are constant, the result follows.
\end{proof}
\end{lemma}

\begin{lemma} \label{S3calc} Let $(S^3, g)$ be the round three-sphere.  Then
given a function $\gp \in C^\infty(S^3)$
\begin{gather}
d \N^* \N \gp = \N^* \N d \gp + 2 d \gp
\end{gather}
and
\begin{gather}
\gd \left( \N^* \N \right) d \gp = - \left(\N^* \N \right)^2 \gp + 2 \N^* \N \gp
\end{gather}
and given a one-form $\ga \in \wedge^1 T^* S^3$
\begin{gather}
\gd \left(\N^* \N \right) \ga = \N^* \N \gd \ga - 2 \gd \ga
\end{gather}
and
\begin{gather}
\gd \left( \N^* \N \right)^2 \ga = \left(\N^* \N \right)^2 \gd \ga - 4 \N^* \N
\gd \ga + 4 \gd \ga
\end{gather}
and
\begin{gather}
\gd L_{\ga^*} g = - \N^* \N \ga + \N \gd \ga + 2 \ga
\end{gather}
and
\begin{gather}
L_{\N^* \N} g = \N^* \N L_\tau g + 4 L_\tau g - 2 \gd \tau g
\end{gather}
and given $B$ a traceless symmetric two-tensor
\begin{gather}
\gd \N^* \N B = \N^* \N \gd B - 4 \gd B + 2 d \tr B
\end{gather}
and
\begin{gather}
\gd \left( \N^* \N \right)^2 = \left(\N^* \N \right)^2 \gd B - 8 \N^* \N \gd B +
16 \gd B
\end{gather}
\begin{proof} We use the convention that $R_{ijij} > 0$, and more specifically
choosing normal coordinates at one point we have
\begin{align*}
R_{ijkl} =&\ \left(g_{ik} g_{jl} - g_{jk} g_{il} \right)\\
R_{ik} =&\ 2 g_{ik}\\
R =&\ 6
\end{align*}
First for a function $\gp$ we compute
\begin{align*}
d \N^* \N \gp =&\ - \N_i \N^k \N_k \gp\\
=&\ - \N^k \N_i \N_k \gp - R_{ikkl} \N^l \gp\\
=&\ - \N^k \N_k \N_i \gp + 2 \N_i \gp\\
=&\ \N^* \N d \gp + 2 d \gp
\end{align*}
Next
\begin{align*}
\gd \left( \N^* \N \right) d \gp =&\ - \N^i \N^j \N_j \N_i \gp\\
=&\ - \N^i \N^j \N_i \N_j \gp\\
=&\ - \N^i \left(\N_i \N^j \N_j \gp + R_{jijk} \N_k \gp \right)\\
=&\ - \left(\N^* \N \right)^2 \gp - R_{jijk} \N^i \N_k \gp\\
=&\ - \left(\N^* \N \right)^2 \gp + 2 \N^* \N \gp
\end{align*}
In the fourth line we used that the curvature is parallel on $S^3$.
Next we compute
\begin{align*}
\gd \left( \N^* \N \right) \ga =&\ - \N^i \N^j \N_j \ga_i\\
=&\ - \N^j \N_i \N_j \ga_i - R_{ijjk} \N^k \ga_i - R_{ijik} \N_j \ga^k\\
=&\ - \N^j \N_i \N_j \ga_i\\
=&\ - \N^j \left(\N_j \N_i \ga_i + R_{ijik} \ga^k \right)\\
=&\ \N^* \N \gd \ga - R_{ijik} \N^j \ga^k\\
=&\ \N^* \N \gd \ga - 2 \gd \ga
\end{align*}
and again we have used that the curvature is parallel.  Using this computation
we have
\begin{align*}
\gd \left( \N^* \N \right)^2 \ga =&\ \N^* \N \gd \N^* \N \ga - 2 \gd \N^* \N
\ga\\
=&\ \N^* \N \left( \N^* \N \gd \ga - 2 \gd \ga \right) - 2 \left(\N^* \N \gd \ga
- 2 \gd \ga \right)\\
=&\ \left(\N^* \N \right)^2 \gd \ga - 4 \N^* \N \gd \ga + 4 \gd \ga
\end{align*}
Also we compute
\begin{align*}
\gd L_{\ga^*} g =&\ \N^i \left( \N_i \ga_j + \N_j \ga_i \right)\\
=&\ - \N^* \N \ga + \N_j \N^i \ga_i + R_{ijil} \ga^l\\
=&\ - \N^* \N \ga + \N \gd \ga + 2 \ga
\end{align*}
Finally we compute
\begin{align*}
L_{\N^* \N \tau} g =&\ - \N_i \N^k \N_k \tau_j - \N_j \N^k \N_k \tau_i\\
=&\ - \N^k \N_i \N_k \tau_j - R_{ikkl} \N^l \tau_j - R_{ikjl} \N^k \tau^l\\
&\ - \N^k \N_j \N_k \tau_i - R_{jkkl} \N^l \tau_i - R_{jkil} \N^k \tau^l\\
=&\ - \N^k \N^k \N_i \tau_j - 2 R_{ikjl} \N^k \tau^l + 2 \N_i \tau_j\\
&\ - \N^k \N^k \N_j \tau_i - 2 R_{jkil} \N^k \tau^l + 2 \N_j \tau_i\\
=&\ \N^* \N L_\tau g + 4 L_\tau g - 2 \gd \tau g
\end{align*}
Now for $B$ a symmetric two-tensor
\begin{align*}
\gd \N^* \N B =&\ - \N^i \N^j \N_j B_{ik}\\
=&\ - \N^j \N^i \N_j B_{ik} - R_{ijjl} \N^l B_{ik} - R_{ijil} \N^j B^l_k -
R_{ijkl} \N^j B_{i}^l\\
=&\ - \N^j \N^i \N_j B_{ik} - R_{ijkl} \N^j B_{i}^l\\
=&\ - \N^j \N^i \N_j B_{ik} - \gd B + d \tr B\\
=&\ - \N^j \left( \N_j \N^i B_{ij} + R_{ijil} B^l_k + R_{ijkl} B_i^l \right) -
\gd B + d \tr B\\
=&\ \N^* \N \gd B - R_{ijil} \N^j B^l_k - R_{ijkl} \N^j B_{il} - \gd B + d \tr
B\\
=&\ \N^* \N \gd B - 4 \gd B + 2 d \tr B
\end{align*}
Using this computation we have
\begin{align*}
\gd \left( \N^* \N \right)^2 =&\ \N^* \N \gd \N^* \N B - 4 \gd \N^* \N B\\
=&\ \N^* \N \left( \N^* \N \gd B - 4 \gd B \right) - 4 \left( \N^* \N \gd B - 4
\gd B \right)\\
=&\ \left(\N^* \N \right)^2 \gd B - 8 \N^* \N \gd B + 16 \gd B
\end{align*}
\end{proof}
\end{lemma}

\bibliographystyle{hamsplain}

\end{document}